# Digitalizing Railway Operations: An Optimization-Based Train Rescheduling Model for Urban and Interurban Disrupted Networks

**Shayan Bafandkar**
The Grainger College of Engineering
Department of Civil and Environmental Engineering
University of Illinois at Urbana-Champaign, Urbana, IL, USA, 61801
Email: sbafan2@illinois.edu
https://orcid.org/0009-0009-8172-5751

**Yousef Shafahi**
Department of Civil Engineering
Sharif University of Technology, Tehran, Iran, 1458889694
Email: shafahi@sharif.edu
https://orcid.org/0000-0003-4267-4348

**Alireza Eslami**
Department of Civil Engineering
Sharif University of Technology, Tehran, Iran, 1458889694
Email: alireza.eslami98@sharif.edu

**Alireza Yazdiani**
Department of Civil and Environmental Engineering
Cornell University, Ithaca, NY, USA, 14853
Email: ay373@cornell.edu
https://orcid.org/0000-0002-4135-5209

**Corresponding author**
Shayan Bafandkar, sbafan2@illinois.edu

**Abstract**

This study introduces a novel methodology for managing train network disruptions across the entire rail network, leveraging digital tools and methodologies. The approach involves two stages, taking into account possible and practical features such as allowing trains to occupy opposite tracks and considering infrastructure capacity for train stops. In the first stage, important nodes within the train network are identified, considering both a topological feature and passenger demand. Subsequently, the network is aggregated based on these important nodes, employing a digital approach to reduce problem complexity. In the second stage, we develop an Integer Programming model for train rescheduling. We then solve this model using the CPLEX solver to evaluate its efficiency. The first case study applies this methodology to the Iranian railway, which is known as a sparse rail network. The results show minimal deviation from the initial train timetable due to the low frequency of trips in each block. Although the approach successfully addresses the train rescheduling problem for various disruption scenarios on the Iranian railway, the excessive computational time required by the optimization model prompts us to make adjustments. Finally, the second case study demonstrates the implementation of the adjusted model in a busy test network. This adaptation significantly reduces computational time by up to 88%. It can be effectively utilized for disruption management in busy networks, where trains need to receive a secondary timetable promptly when facing disruptions.

## 1. Introduction

Nowadays, railroad systems are widely used for urban trips due to their low cost and relatively high level of service, especially by commuters in developing countries such as Poland [1] and Russia [2]. Moreover, the emergence of rapid railroads has made rail systems even more competitive compared to other transportation modes, such as air travel. The growing demand for rail transportation, along with disruptions caused by natural disasters or malicious attacks, has highlighted the importance of rail network disruption management in recent transportation engineering studies. One approach that offers greater stability and more efficient operation in the rail transportation sector is the digitalization of rail networks ([3], [4]).

There is no universally accepted definition for the digitalization of transportation systems, and this novel concept has recently been studied from various perspectives. However, most available definitions include the introduction of smart approaches into transportation network management, with an emphasis on automation and intelligence. In the case of rail networks, Poliński and Ochociński [1] divide previous digitalization efforts into five major trends: connected commuters, mobility as a service (MAAS), project management as a service (PMAAS), automation and integration of rail traffic management systems, and the internet of trains. Our study falls within the last area, as it considers rail traffic safety and resilience improvement and focuses on the system's ability to absorb the negative effects of disruptions and rapidly achieve a new equilibrium status with the aid of digital approaches and an optimization model. Since one of the main goals of rail network digitalization is to enhance the efficiency of railway operations and reduce delays, we formulated the objective function of the optimization model introduced in this study to minimize the total travel time of the trains, which, in turn, reduces delays caused by disruptions in the rail network.

To be more specific, this study focuses on presenting an optimization model for solving the train rescheduling problem in a disrupted railway network, with the aim of facilitating network disruption management. On one hand, disruption management requires a prompt response to be effective on a network-wide scale. On the other hand, solving the train rescheduling problem is a highly complex task when considering the entire train network. To address these challenges, we follow a digital and smart approach consisting of two steps. First, important nodes within the network are identified using nodes' degree and passenger demand. These nodes play a crucial role in network stability, and any disturbance to them can significantly impact overall network performance. Given the complexity and vastness of railway networks, we employ a smart and computationally efficient approach to summarize them into average-scale networks by isolating these important nodes. Second, we apply the rescheduling optimization model to the average-scale network to provide a solution. This process highlights the use of optimization tools and methodologies within a digital transportation system.

Notably, our approach enables rescheduling optimization across the entire network, unlike conventional methods that focus on individual blocks. Additionally, we incorporate practical and feasible features into the model formulation, such as allowing trains to occupy opposite tracks and accounting for infrastructure capacity at train stops. These considerations enhance the applicability of our results to real-world scenarios. We validate the effectiveness of our methodology through two case studies. The first case study focuses on the Iranian Rail Network, a sparse rail network, and demonstrates minimal deviation from the initial train timetable due to the low frequency of trips in each block. In the second case study, we apply the adjusted model to a crowded test network, achieving a significant reduction in computational time by up to 88%. This adaptation is particularly valuable for disruption management in busy networks, where trains must promptly receive a secondary timetable during disruptions.

The following is divided into six sections. In Section 2, previous related studies are briefly reviewed. The optimization model and the solution approach are introduced in Sections 3 and 4, respectively. Two case studies are presented in Sections 5 and 6. An adjustment to the optimization model to improve its efficiency is discussed in Section 7. Finally, conclusions are provided in Section 8.

## 2. Literature Review

Jing et al. [5], Wang et al. [6], and Sun et al. [7] all proposed methods for identifying critical nodes in a railway system. Jing et al. [5] used the mean excess criticality problem (MECP) metric, which considers both the network's topology and passenger flow, while Wang et al. [6] utilized complex transportation network theory to measure criticality. Sun et al. [7] used a topological and behavioral approach to determine the critical nodes of Beijing's metro network by gathering passenger information using the Automated Fare Collection (AFC) system.

Furthermore, several studies focused on optimizing train rescheduling, such as Bešinović et al. [8], proposed the Infrastructure Restoration and Transport Management (IRTM) model, which measures rail network resilience using an optimization approach and divides the original network into four subnetworks: infrastructure network, railway network, passenger network, and restoration network. Tang et al. [9] studied rail network resilience from a quantitative point of view and expanded on Xu and Ng [10]'s linear model by taking commuters' transfer walking time into account. Yang et al. [11] proposed a two-stage fuzzy integer optimization model that determines the optimal sequence of trains' movement, and Gong et al. [12] formulated an integer nonlinear model for train timetabling. Bärmann et al. [13] proposed a train timetabling model to reduce energy consumption and achieved significant reductions in railway instances. Sessa et al. [14] presented a hybrid stochastic approach for creating train trajectory, ensuring its robustness. Wang et al. [15]

evaluated rail systems' safety using the cusp catastrophe model. Leng et al. [16] compared two methodologies used for the rescheduling process of passengers, train operators, and infrastructure managers in railway disruptions. Herrigel et al. [17] modeled the railway timetabling problem as a periodic event scheduling problem (PESP) and solved this model using a novel methodology called sequential decomposition. Toletti et al. [18] presented a model for railway rescheduling based on the resource conflict graph model introduced by Caimi [19]. Ghaemi et al. [20] classified the current railway disruption management practices performed in the Netherlands with regard to the 3-phased bathtub model and then showed the effectiveness of a microscopic rescheduling model when applied to a corridor of the Dutch railway. Kecman et al. [21] mentioned four macroscopic rescheduling models to modify the rail network's operations during disruption, relying on the timed event graphs previously introduced in Goverde [22] & [23].

Bešinović [24] conducted a comprehensive review of available literature related to railway resilience. Ge et al. [25] provided a survey on the different concepts of disruption and the possible approaches to manage them by expounding on the effects of a disruption on the performance of public transport networks, including rail transits, over time. Since the main focus of this study is to introduce a mathematical programming model for the train rescheduling problem, previous works of research on the railway resilience are not extensively discussed.

To the best of our knowledge, none of the past studies has applied the concepts of important node identification, mixed single-track and double-track operations, and infrastructure capacity simultaneously to solve the train scheduling problem considering the whole network.

## 3. Problem Setting and formulation

Train timetables regulate the flow of trains in the network based on details such as the initiation time of train movements, stations at which each train should stop, and travel time for blocks traversed by each train. Proper management of train timetables is crucial to satisfy safety standards and transportation demands of travelers and commodities. Otherwise, problems such as delay propagation, spoilage of perishable products, and traveler dissatisfaction may occur, resulting in significant financial losses for railroad beneficiaries. Since railroads play a vital role in a country's economy, they should run based on an optimal schedule during normal conditions and a proper model for rescheduling during disruptions to ensure stable performance.

In this study, we formalize an optimization model for the train rescheduling problem, focusing on the train dispatch management, that considers the entire network, safety issues, infrastructure capacity, and block disruption. We explain the details of this model thoroughly below.

### 3.1 Problem Definition and Assumptions

Our model solves the trains rescheduling problem by minimizing the total travel time of trains in the disrupted network according to an initial timetable while considering the mixed single and double-track operations, making it more useful for real-world networks. The model makes use of the space-time trajectory concept, introduced by Yang et al. [11], where stations are nodes and blocks are arcs, allowing it to divide the time horizon into several intervals of equal length and take into account the physical features of the subject network simultaneously. An example of the space-time trajectory is depicted in Figure 1.

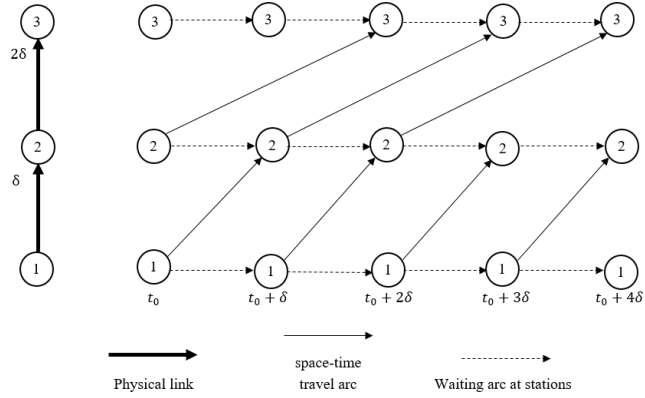

**Figure 1-** An example of the space-time trajectory introduced by Yang et al. [11]

The following assumptions are considered in this study:

- Disruptions occur in blocks, not nodes, and only complete disruptions are taken into account.
- Trains are divided into two groups - positive or negative - based on their direction, allowing the model to consider mixed single and double-track operations.
- Trains can only overtake each other while at stations; therefore railway switches are available in all of the stations.
- Each block can be occupied by several trains at the same time if the minimum safe headway is available.
- Two trains can only conflict in double blocks or at stations.
- An initial and satisfactory timetable for an undisrupted condition is available.
- Minimum dwelling time should be satisfied at each station.
- All trains move at a constant speed that is the same among all of them.
- All trains present in the network during the time horizon are assumed to be passenger trains. This would make the model applicable to urban railways as well as inter-city networks.

### 3.2 Decision Variables and Parameters

The model's decision variables are two types of binary variables:

- $x_{ij}^{k,tt'}$ is a binary variable indicating the occupation status of block ij. If block ij is occupied by train k in the tt' time period, $x_{ij}^{k,tt'} = 1$; otherwise, $x_{ij}^{k,tt'} = 0$.
- $p_{ij}^{k,m}$ is a binary variable that is set to 1 if train k has priority to enter block ij over train m; otherwise, it is equal to zero.

Additionally, the following four auxiliary variables are defined to make the model's constraints more comprehensible:

- The times when train k enters and departs block ij are represented by $a_{ij}^k$ and $d_{ij}^k$, respectively.

• The times when train k enters and departs node i are represented by $a_i^k$ and $d_i^k$, respectively. Other model's parameters are discussed in Table 1.

**Table1**
Parameters used in the rescheduling optimization model:

| notation | Definition |
|---|---|
| $N$ | Set of network nodes |
| $A$ | Set of network blocks |
| $A^*$ | Set of disrupted blocks |
| $R^+, R^-$ | Set of north and south headed trains |
| $R^*$ | Set of trains directly affected by the disruption |
| $RU^k$ | Set of blocks that should be traversed by train k |
| $O^k, D^k$ | Origin and destination nodes of train k |
| $h_{ij}$ | Minimum headway between two consecutive trains on block ij. |
| $DT_i, Cap_{ij}$ | Minimum dwelling time and capacity at block ij |
| $\delta, T$ | Time interval length and horizon under review |
| $t_0, t_{dis}, \tau_{dis}$ | Commencement time of subject horizon, disruption start time, and disruption duration |
| $\bar{a}_{ij}^k, \bar{d}_{ij}^k, \bar{x}_{ij}^k$ | Actual time when train k enters/leaves block ij in initial timetable, and binary parameter indicating occupation status of block ij in initial timetable |
| $\bar{a}_i^k, \bar{d}_i^k$ | Actual time when train k enters/departs node i in initial timetable |
| $dis_{ij}^{tt'}$ | An integer parameter indicating the status of block ij according to the number of its tracks disrupted over the tt' time period. In this parameter t' equals t+$\tau_{dis}$ |
| $bt_{ij}$ | Block type binary indicator which is equal to 1 if block ij has two tracks or 0 if it has only one track |
| $\beta, Z, M$ | Delay threshold, total number of time intervals which equals the time horizon ($T$) divided by the time interval length ($\delta$), and a very big number |

Eventually, the objective function and all of its constraints are mentioned in Table 2. The description of each equation is provided.

**Table2**
Objective function and constraints of the rescheduling optimization model:

| Equation | No. |
|---|---|
| $Min \sum_k \left( a_{D^k}^k - a_{O^k}^k \right)$ | 1 |
| $\sum_{ij \in RU^k} x_{ij}^{k,tt'} \leq 1; \ \forall t, t' \in T; \quad \forall k \in \{R^+ \cup R^-\}$ | 2 |
| $\sum_{t'} \sum_t x_{ij}^{k,tt'} = 0; \ \forall k \in \{R^+ \cup R^-\} \quad \forall ij \in A - RU^k, i \neq j$ | 3 |
| $\sum_{ij \in RU^k} x_{ij}^{k,tt'} = 0; \ \forall k \in \{R^+ \cup R^-\} \quad \forall t, t' \in T; \ t' < \bar{a}_{O^k}^k \ or \ t > t_0 + Z \times \delta$ | 4 |

$$\sum_{ji'} x_{ji'}^{k,t't''} - \sum_{ij} x_{ij}^{k,tt'} = \begin{cases} 1; if\ j=O^k,\ t'=\bar{d}_{O}^{k} \\ -1; if\ j=D^k \\ 0; otherwise \end{cases} ; \forall k \in \{R^+ \cup R^-\} \quad \forall ij, ji' \in A \quad \forall t, t' \in T \qquad 5$$

$$x_{ij}^{k^+,tt'} + x_{ij}^{k^-,ll'} \leq 1 + bt_{ij} - dis_{ij}^{gg'}$$
$$\forall t, t', l, l', g, g' \in T, \ tt' \cap ll' \cap gg' \neq \emptyset \quad \forall ij \in RU^{k^+} \cap RU^{k^-} \quad \forall k^+ \in R^+ \quad \forall k^- \in R^- \qquad 6$$

$$\sum_{g} \sum_{g'} x_{ij}^{k,gg'} \times (g'-g) = \sum_{t} \sum_{t'} \bar{x}_{ij}^{k,tt'} \times (t'-t)$$
$$\forall t, t', g, g' \in T \quad \forall ij \in A, i \neq j \quad \forall k \in \{R^+ \cup R^-\} \qquad 7$$

$$a_{ij}^k = \sum_{t} \sum_{t'} x_{ij}^{k,tt'} \times t; \ \forall k \in \{R^+ \cup R^-\} \quad \forall ij \in RU^k \qquad 8$$

$$d_{ij}^k = \sum_{t'} \sum_{t} x_{ij}^{k,tt'} \times t'; \ \forall k \in \{R^+ \cup R^-\} \quad \forall ij \in RU^k \qquad 9$$

$$a_{j}^k = \sum_{t'} \sum_{t} x_{ij}^{k,tt'} \times t'; \ \forall k \in \{R^+ \cup R^-\} \quad \forall ij \in RU^k \quad \forall j \in N^k (j \neq O^k) \qquad 10$$

$$a_{O^k}^k = \bar{a}_{O^k}^k; \ \forall k \in \{R^+ \cup R^-\} \qquad 11$$

$$d_{j}^k = \sum_{t} \sum_{t'} x_{ij}^{k,tt'} \times t' + \sum_{t} \sum_{t'} x_{jj}^{k,tt'} \times (t'-t); \ \forall k \in \{R^+ \cup R^-\} \quad \forall ij \in RU^k \quad \forall i \in N^k \qquad 12$$

$$p_{ij}^{k,m} + p_{ij}^{m,k} = 1; \ \forall ij \in \{RU^k \cap RU^m\} \quad \forall k, m \in \{R^+ \cup R^-\}, k \neq m \qquad 13$$

$$a_{ij}^k + h_{ij} \leq a_{ij}^m + M \times (1 - p_{ij}^{k,m}); \ \forall ij \in \{RU^k \cap RU^m\} \quad \forall k, m \in \{R^+ \cup R^-\} \qquad 14$$

$$a_{ij}^m + h_{ij} \leq a_{ij}^k + M \times (1 - p_{ij}^{m,k}); \ \forall ij \in \{RU^k \cap RU^m\} \quad \forall k, m \in \{R^+ \cup R^-\} \qquad 15$$

$$d_{ij}^k + h_{ij} \leq d_{ij}^m + M \times (1 - p_{ij}^{k,m}); \ \forall ij \in \{RU^k \cap RU^m\} \quad \forall k, m \in \{R^+ \cup R^-\} \qquad 16$$

$$d_{ij}^m + h_{ij} \leq d_{ij}^k + M \times (1 - p_{ij}^{m,k}); \ \forall ij \in \{RU^k \cap RU^m\} \quad \forall k, m \in \{R^+ \cup R^-\} \qquad 17$$

$$d_{i}^k - a_{i}^k \geq DT_i; \ \forall i \in N^k \quad \forall k \in \{R^+ \cup R^-\} \qquad 18$$

$$\sum_{k} x_{ii}^{k,tt'} + \sum_{k} x_{jj}^{k,tt'} \leq Cap_{ij}; \ \forall ij \in RU^k \quad \forall t, t' \in T \qquad 19^{[1]}$$

$$\sum_{i} \left( d_{i}^k - \bar{d}_{i}^k \right) \leq \beta; \ \forall i \in N^k \quad \forall \bar{d}_{i}^k \geq t_{dis} \quad \forall k \in \{R^+ \cup R^-\} \qquad 20$$

$$d_{i}^k - \bar{d}_{i}^k = 0; \ \forall i \in N^k \quad \forall \bar{d}_{i}^k < t_{dis} - \beta \quad \forall k \in \{R^+ \cup R^-\} \qquad 21$$

$$x_{ij}^{k,tt'} - \bar{x}_{ij}^{k,tt'} = 0; \ \forall ij \in A \quad \forall t, t' < t_{dis} - \beta \quad \forall k \in \{R^+ \cup R^-\} \qquad 22$$

$$d_{i}^k - \bar{d}_{i}^k \geq 0; \ \forall i = N^k \quad \forall k \in \{R^+ \cup R^-\} \qquad 23$$

$$a_{i}^k - \bar{a}_{i}^k \geq 0; \ \forall i = N^k \quad \forall k \in \{R^+ \cup R^-\} \qquad 24$$

---

[1] The capacity is defined in terms of the block's capacity to offset the impact of removing non-critical nodes that may be available within each block. It is worth mentioning that this constraint could be replaced by "$\sum_k x_{ii}^{k,tt'} \leq Cap_i$" to consider node capacity.

The model's objective function is formulated based on the total travel time for all trains during the time horizon and represented in Equation 1.

The model's constraints can be classified into ten classes:

(1) *Constraints Related to the Viability of the Rescheduled Timetable:* The rescheduled timetable should provide a feasible plan in several aspects. Constraint 2 ensures that each train cannot occupy more than one block in the space-time trajectory at the same time. Constraint 3 confines each train's route to its original one introduced in the given initial timetable. Constraint 4 prevents trains from continuing their path after the time horizon ($T$) has been reached or commencing their trips before the initiation of the initial timetable.

(2) *Constraints Related to Flow Balance in Blocks:* Constraint 5 provides a continuous and balanced train flow in the network.

(3) *Constraints Related to Trains' Conflict*: Constraint 6 prevents trains with conflicting directions from occupying a single block simultaneously. This constraint enables the model to consider two-way double blocks and also to provide a chance for trains affected by the disruption to use their opposing track if possible. In other words, if a double-track block faces a disruption in which one of its tracks becomes disrupted, it could be considered as a single-track block in the disruption period and trains with opposing directions could still use that block as long as their paths do not conflict.

(4) *Constraints Related to Blocks' Travel Time:* Each block's travel time is taken into account in constraint 7. This constraint is written assuming that all trains travel at a constant speed.

(5) *Constraints Related to Blocks' Auxiliary Variables:* The actual time when train k enters and leaves block ij is defined in constraints 8 and 9, respectively. These two variables will be used in writing the blocks' headway constraints.

(6) *Constraints Related to the Modes' Auxiliary Variables: The* actual time when train k enters and leaves node i is defined in Constraints 10 and 12, respectively. Moreover, Constraint 11 makes each train's arrival time at its origin node match the arrival time defined in the given initial timetable.

(7) *Constraints Related to Trains' Succession in Blocks:* The order of trains' movement in different blocks is defined in Constraint 13, then Constraints 14 to 17 are defined to provide a safe headway between successive trains in each block.

(8) *Constraints Related to Dwelling Time at Each Node:* Passenger trains should stop at predetermined nodes to embark (disembark) passengers and cargo, and also undergo safety checks. Constraint 18 ensures that trains' stoppage time is more than a minimum interval called dwelling time which may have distinctive values for different nodes. Constraint 19 takes into account each node's dwelling capacity.

(9) *Constraints Related to Divergence from the Given Initial Timetable:* Constraint 20 limits each train's delay in the rescheduled plan to a reasonable threshold according to its delay in the given initial timetable and a coefficient (β) that should be determined based on the potential disruption's length. This particular constraint simplifies the calculation process since it limits the feasible area of the problem. Moreover, Constraints 21 and 22 ensure that trains' movements, before the commencement of the disruption, accord with those determined in the given initial timetable.

(10) *Constraints Related to Optimality of the Given Initial Timetable*: Constraints 23 and 24 indicate that the rescheduled plan's total delay cannot be smaller than that of the given initial timetable.

## 4. Model's solution approach

The introduced formulation is an integer programming (IP) model that could be readily solved by commercial optimization software such as GAMS and LINGO. This study makes use of the CPLEX solver that relies on several algorithms including enumeration techniques such as the branch and bound approach as well as cutting-plane techniques to find the optimum solution. Although this solver can accurately determine the optimal solution of linear programming models, it has several drawbacks related to the amount of memory and time it needs. Therefore, it cannot solve large-scale problems such as the train rescheduling problem, and an algorithm is required to summarize these large-scale problems into an average-scale network by considering only the rail network's critical nodes. Such an algorithm is introduced in Section ۴/۱ to make the mentioned model solvable for the CPLEX solver.

### 4.1 Model' solution approach

Railroads are comprised of many nodes, some of which are created solely for political and social preferences. Hence, their performance does not affect the whole network's level of service. On the other hand, some nodes are of utmost importance because of their high demand or significant economic impacts, making them crucial in determining the whole network's performance. The algorithm introduced in this study, depicted in Figure 2, pinpoints the important nodes of the subject rail network by considering topological factors as well as passenger demand. The algorithm calculates a criticality index, consisting of two terms indicated in Equations 26 and 27, for each node. Equation 26 measures each node's normalized degree, while Equation 27 calculates each node's normalized passenger demand.

$$\bar{D}_i = \begin{cases} \dfrac{d_i - d_{min}}{d_{max} - d_{min}} \text{; if } d_i \neq d_{min} \\ \dfrac{0.5}{d_{max} - d_{min}} \text{; if } d_i = d_{min} \end{cases} \tag{26}$$

$$\bar{P}_i = \frac{p_i - 0.8 \times p_{min}}{p_{max} - p_{min}} \tag{27}$$

Where $\bar{D}_i$ and $d_i$ are the normalized measure of the node i's degree and degree of node i, respectively. $d_{min}$ and $d_{max}$ indicate the minimum and maximum degree among the network's nodes, respectively. Moreover, $\bar{P}_i$ and $p_i$ are indicators of the normalized measure of the node i's passenger demand and its monthly passenger demand. The minimum and maximum passenger demand among all the network's nodes are considered as $p_{min}$ and $p_{max}$.

Eventually, the importance index of each node is calculated according to Equation 28, in which the mentioned normalized measures are combined based on two importance power weights ($\alpha_1$, $\alpha_2$). These power weights have values between zero and one, and their summation should be equal to one.

$$Cr_i = (\bar{D}_i)^{\alpha_1} \times (\bar{P}_i)^{\alpha_2} \times 100 \tag{28}$$

In the end, the large-scale network can be aggregated to a medium-scale network consisting of nodes with the highest importance index, connected with a two-way double track or a two-way single track based on their status in the large-scale network. Moreover, trains originating from omitted nodes will not be considered in the solution process. Power weights ($\alpha 1$, $\alpha 2$) are determined based on the network's structure and decision-makers' goals. The higher these power weights are, the lower the importance of their corresponding measures.

Although omitting the non-important nodes in the train rescheduling process would negatively impact the model's results as these nodes may provide additional capacity for meet-pass and overtaking maneuvers, the benefit of omitting them, which is severe reduction in computational time, would offset this drawback as these nodes usually have low passenger demand with low trip frequency.

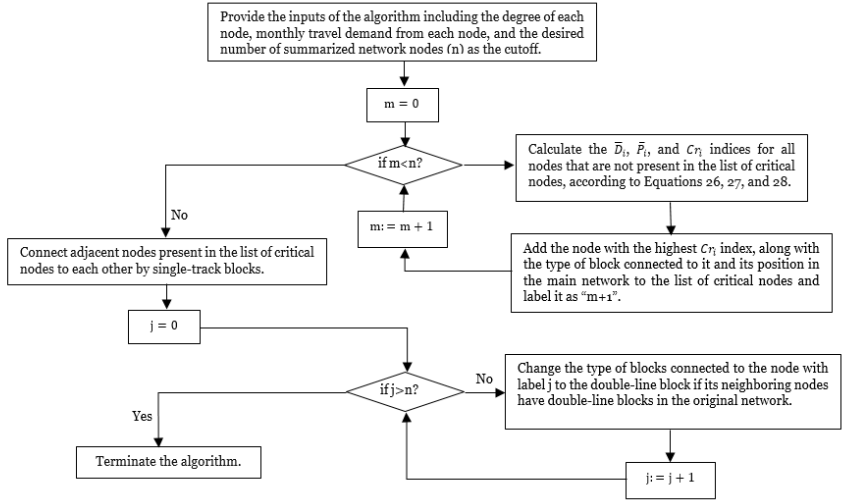

**Figure 2-** The flowchart of the important nodes algorithm

## 5. Case Study 1: Iran's Railway

In this section, the optimization model and algorithm introduced in Sections 3 and 4 have been implemented on Iran's railway Network to validate their competency. As of October 2022, Iran's railroad includes 12,652 km of single blocks and 1,426 km of double blocks; therefore, it is not considered a dense network.

### 5.1 Implementation of the Critical Nodes Algorithm

In railroads with a low level of density, passenger flow has the utmost importance in determining important nodes because the node degree for most of the nodes is more or less equal to 2. Therefore, $\alpha 1$ and $\alpha 2$ are considered 0.6 and 0.4, respectively. The cutoff was considered nine nodes because this sufficed to create an aggregated network covering a vast area of the rail network in the country. In other words, the cutoff creates an aggregated network that spans in all the geographical directions (North, South, East, West) and it can resemble the whole Iranian network. Moreover, the timetable of the trains of Iran's railway in October 2022 was obtained from the official website of Iran's railway (https://www.rai.ir/). Then, nine important nodes of that network were determined by the important nodes algorithm and written in bold in Table 3. The mentioned important nodes are so crucial that they can solely determine the level of service of the rail network; hence they can be used to create a medium-scale network. Such a network is demonstrated in Figures 3 and 4.

**Table3**

Iran's railway important nodes:

| Station's Name | Passengers monthly demand | Node's degree | $\bar{D}_i$ | $\bar{P}_i$ | $\overline{Cr}_i$ |
|---|---|---|---|---|---|
| **Tehran** | **2556200** | **3** | **1** | **0.738185432** | **88.6** |
| Khaf | 44800 | 2 | 0.5 | 0.004115515 | 7.3 |
| Tabas | 44800 | 2 | 0.5 | 0.004115515 | 7.3 |
| **Mashhad** | **3459600** | **1** | **0.25** | **1.002244826** | **43.6** |
| **Isfahan** | **230400** | **2** | **0.5** | **0.058365486** | **21.2** |
| Shiraz | 153600 | 1 | 0.25 | 0.035917222 | 11.5 |
| **Tabriz** | **499200** | **2** | **0.5** | **0.136934409** | **29.8** |
| Urmia | 96000 | 1 | 0.25 | 0.019081024 | 8.9 |
| **Ahvaz** | **384000** | **3** | **1** | **0.103262013** | **40.3** |
| Zanjan | 136200 | 2 | 0.5 | 0.030831287 | 16.4 |
| Ghazvin | 57600 | 2 | 0.5 | 0.007856892 | 9.5 |
| **Karaj** | **172800** | **2** | **0.5** | **0.041529288** | **18.5** |
| Kashan | 96000 | 2 | 0.5 | 0.019081024 | 13.5 |
| **Qom** | **211200** | **3** | **1** | **0.05275342** | **30.8** |
| Arak-Qom | 38400 | 2 | 0.5 | 0.002244826 | 5.8 |
| Malayer | 96000 | 2 | 0.5 | 0.019081024 | 13.5 |
| Kermanshah | 96000 | 1 | 0.25 | 0.019081024 | 8.9 |
| Yazd | 144900 | 2 | 0.5 | 0.033374255 | 16.9 |
| Bandar Abbas | 38400 | 1 | 0.25 | 0.002244826 | 3.8 |
| **Kerman** | **230400** | **2** | **0.5** | **0.058365486** | **21.2** |
| Sari | 96000 | 2 | 0.5 | 0.019081024 | 13.5 |
| Hamedan | 134400 | 1 | 0.25 | 0.030305156 | 10.7 |
| Rasht | 164100 | 1 | 0.25 | 0.038986321 | 11.9 |
| Khoramshahr | 76800 | 1 | 0.25 | 0.013468958 | 7.8 |
| **Maraghe** | **96000** | **3** | **1** | **0.019081024** | **20.5** |
| Mianeh | 42000 | 3 | 1 | 0.003297089 | 10.2 |

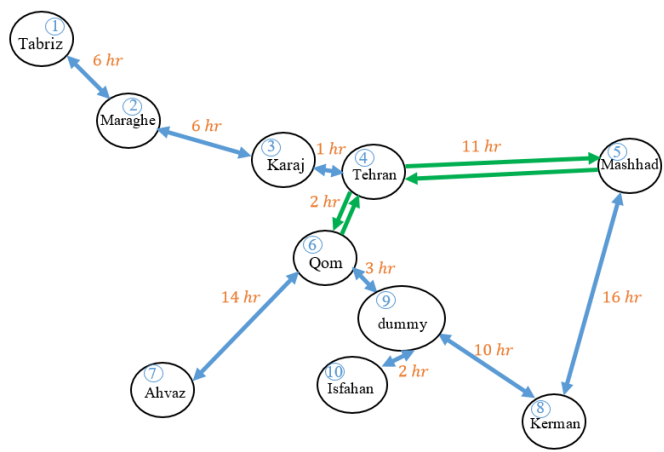

**Figure 3-** The aggregated network of Iran's railway

It is worth noting that a dummy node is added to the aggregated network because of the overlap of Qom–Isfahan and Qom–Kerman trains which should be separated at a point so that their trips can be divided at that point. The dummy node can serve this purpose and is located in Kashan city on the map; hence, it could also be labeled as "Kashan".

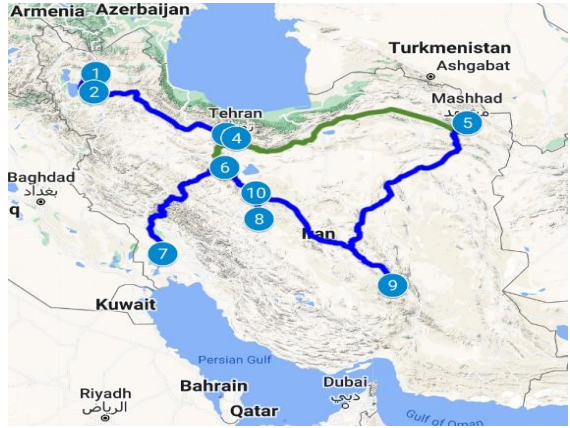

**Figure 4-** The exact configuration of the aggregated network captured by Google Maps 2023

### 5.2 Implementation of the Optimization Model for Train Rescheduling

The network portrayed in Figure 3 is comprised of nine important nodes and a dummy node that separates Qom–Kerman and Qom–Isfahan trains. Furthermore, these nodes are connected by blue and green arrows indicating single blocks and double blocks, respectively. Each block's travel time is also written above its corresponding arrow. It is assumed that all trains stop only at their destinations after departing their origin nodes, except the ones that cross the dummy node within their routes. The mentioned trains should stop at the dummy node for 10 minutes regardless of their destinations. Moreover, the minimum safe headway is considered 5 minutes in all blocks. The parameter related to the time horizon (T) is assumed to be 45 hours to cover all trips in the initial timetable, and it starts at 0:00 AM of an even day and ends at 9 PM of the following day. The delay threshold is equal to 480 min, which is logical since such a value can ensure the feasibility of the problem and also prevent passengers' discontent and trip cancellation in Iranian railway. The basic scenario of disruption is a total disruption in Block 46[1] (Tehran–Qom block), because around 40% of the trips should pass this block. The disruption period starts at hour 12 and ends at hour 14 of the first day under review ($dis_{46}^{720,840} = 2$). In the end, the initial timetable of trains is presented in Table4.

**Table4**

The initial timetable of trains in Iranian railway network:

| Train No. | Departure Time | Arrival Time | Trip Length (hr) | Departure Frequency | Origin | Destination |
|---|---|---|---|---|---|---|
| 1 | 1:00 AM | 2:00 PM | 13 | Even Days | Mashhad | Qom |

---

[1] The notation used to indicate each block in the first case study is the number of its ending nodes, according to Figure 3.

| 2 | 6:20 AM | 6:20 PM | 12 | Even Days | Mashhad | Karaj |
| 3 (3') | 7:40 AM | 6:40 PM (+1 day) | 11 | Every Day | Mashhad | Tehran |
| 4 (4') | 8:00 AM | 7:00 PM (+1 day) | 11 | Every Day | Tehran | Mashhad |
| 5 | 8:15 AM | 2:25 AM | 18 | Even Days | Isfahan | Mashhad |
| 6 | 10:50 AM | 11:50 PM | 13 | Even Days | Mashhad | Qom |
| 7 | 1:25 PM | 4:35 AM | 15 | Even Days | Kerman | Tehran |
| 8 | 2:00 PM | 2:00 PM | 24 | Every Day | Tabriz | Mashhad |
| 9 | 2:25 PM | 8:35 AM | 18 | Even Days | Mashhad | Isfahan |
| 10 | 3:10 PM | 7:10 AM | 16 | Even Days | Ahvaz | Tehran |
| 11 | 4:10 PM | 4:10 AM | 12 | Even Days | Karaj | Mashhad |
| 12 | 4:55 PM | 5:55 AM | 13 | Every Day | Tabriz | Tehran |
| 13 | 5:25 PM | 9:25 AM | 16 | Even Days | Kerman | Mashhad |
| 14 | 5:45 PM | 12:45 AM | 7 | Even Days | Maraghe | Tehran |

Figures 5 and 6 demonstrate the space-time trajectory for the trains' initial timetable in the summarized undisrupted and disrupted network, respectively. The disruption period is indicated by an orange hatch in Figure 5. Please note that in all of the following space-time diagrams, passing a station is indicated only by vertices.

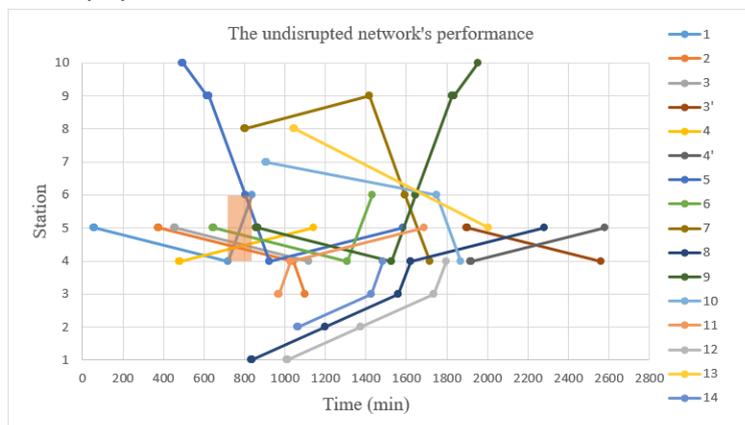

**Figure 5-** The space-time trajectory for the trains' initial timetable in the summarized undisrupted network

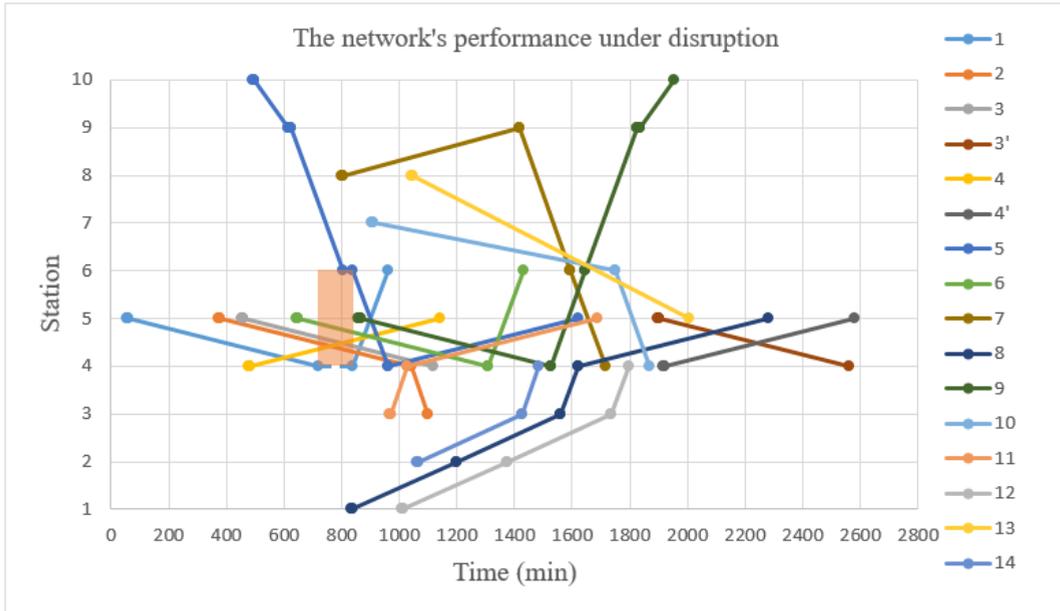

**Figure 6-** The space-time trajectory for the trains' timetable under disruption

It is noteworthy that a total disruption of Block 46 (Tehran–Qom) directly affects the trip schedule of trains 1 and 5, which in turn deteriorates the whole network's performance by increasing the objective function from 13375 minutes to 13530 minutes. Therefore, it can be concluded that a 120-minute disruption interval can cause 155 minutes of delay. The model was solved by the CPLEX solver on a personal computer with a 16 GB RAM, and 2.6 GHz Core i7 CPU in 2386 seconds.

### 5.3 Different Disruption Scenarios

In this section, different disruption scenarios are introduced, for each of which the optimization model is implemented. The model's results for these scenarios, including Blocks 46, 34, 45, and 12, are mentioned in Table 5. These results indicate that occupying the opposite tacks can enhance the network's performance significantly by comparing the total delay in the full and partial disruption of Blocks 46 and 45. For instance, the partial disruption of Block 46 (Tehran to Qom) did not change the initial timetable at all because trains 1 and 5 that were supposed to be affected could use their opposing blocks to complete their travels. Moreover, the total disruption of Block 45 (Tehran–Mashhad) can cause the most severe delay; hence it can be named the most critical block of the network.

**Table5**
The model's results for different disruption scenarios:

| Disrupted block | Disruption period (minute) | The optimization model's objective function | Total delay (minute) | Computational time (second) |
|---|---|---|---|---|
| Block 46's Tehran to Qom link | 720-840 | 13410 | 35 | 2220 |

| Block 45's Mashhad to Tehran link | 720-840 | 14420 | 1045 | 2473 |
|---|---|---|---|---|
| Block 45 | 720-840 | 14780 | 1405 | 2514 |
| Block 34 | 1450-1570 | 13535 | 160 | 1962 |
| Block 12 | 910-1030 | 13585 | 210 | 2185 |
| Block 23 | 1300-1420 | 14010 | 635 | 1909 |

### 5.4 Sensitivity Analysis

One of the introduced optimization model's parameters that ensured a safe and feasible timetable is its minimum headway. This study assumed a constant headway in all blocks, whereas different blocks of a real-world railway can have variable minimum headways. Hence, a sensitivity analysis should be conducted to ensure that the model's results were not significantly affected by the aforementioned assumption. The results of the model for a 10-minute headway and a 20-minute headway are presented in Tables 6 and 7, respectively. The only impact of increasing the minimum headway from 5 to 10 minutes is that train 9 would be affected by the disruption of Block 45. Increasing the minimum headway to 20 minutes does not result in significant changes in the final outcomes.

Another assumed parameter that requires analysis is the duration of disruption, as natural and human-made disruptions are highly unpredictable. As a result, two additional scenarios involving 60-minute and 180-minute disruptions for each of the previously disrupted blocks have been tested, and the findings are shown in Tables 8 and 9. The disruption period in each scenario is selected so that the disruption happens at the busiest time interval of the corresponding block. The results indicate that altering the disruption period has a negligible effect in general. However, its impact on scenarios where the disrupted block has a travel time higher than the disruption period, such as Block 45, is more tangible.

**Table 6**

The model's results for $h_{ij} = 10$ minutes:

| Disrupted block | Disruption period (minute) | $h_{ij}$ (minute) | The optimization model's objective function | Total delay (minute) | Computational time (second) |
|---|---|---|---|---|---|
| Block 46 | 720–840 | 10 | 13530 | 155 | 2431 |
| Block 46's Tehran to Qom link | 720–840 | 10 | 13410 | 35 | 2265 |
| Block 45's Mashhad to Tehran link | 720–840 | 10 | 14440 | 1065 | 2561 |
| Block 45 | 720–840 | 10 | 14800 | 1425 | 2693 |
| Block 34 | 1450–1570 | 10 | 13540 | 165 | 2053 |
| Block 12 | 910–1030 | 10 | 13590 | 215 | 2217 |
| Block 23 | 1300–1420 | 10 | 14025 | 650 | 2038 |

**Table7**
The model's results for $h_{ij} = 20$ minutes:

| Disrupted block | Disruption period (minute) | $h_{ij}$ (minute) | The optimization model's objective function | Total delay (minute) | Computational time (second) |
|---|---|---|---|---|---|
| Block 46 | 720-840 | 20 | 13530 | 155 | 2586 |
| Block 46's Tehran to Qom link | 720-840 | 20 | 13410 | 35 | 2298 |
| Block 45's Mashhad to Tehran link | 720-840 | 20 | 14470 | 1095 | 2701 |
| Block 45 | 720-840 | 20 | 14830 | 1455 | 2782 |
| Block 34 | 1450-1570 | 20 | 13550 | 175 | 21۳۰ |
| Block 12 | 910-1030 | 20 | 13600 | 225 | 2199 |
| Block 23 | 1300-1420 | 20 | 14055 | 680 | 2074 |

As could be inferred from the above tables, altering the headway has not made a significant change in the total delay. This is justifiable because the Iranian railway is a sparse network where the trip frequency at each block is low. As a result, altering the minimum safe headway does not significantly impact the network's performance.

**Table8**
The model's results for a 60-minute disruption:

| Disrupted block | Disruption period (minute) | $h_{ij}$ (minute) | The optimization model's objective function | Total delay (minute) | Computational time (second) |
|---|---|---|---|---|---|
| Block 46 | 770-830 | 5 | 13510 | 135 | ۱۸74 |
| Block 46's Tehran to Qom link | 770-830 | 5 | 13400 | 15 | ۱۷89 |
| Block 45's Mashhad to Tehran link | 750-810 | 5 | 14375 | 1000 | ۱۹44 |
| Block 45 | 750-810 | 5 | 14715 | 1340 | 2۱03 |
| Block 34 | 1510-1570 | 5 | 13390 | 15 | 1۵56 |
| Block 12 | 940-1000 | 5 | 13525 | 150 | 1۶34 |
| Block 23 | 1330-1390 | 5 | 13920 | 545 | 1۴99 |

**Table9**
The model's results for a 180-minute disruption:

| Disrupted block | Disruption period (minute) | $h_{ij}$ (minute) | The optimization model's objective function | Total delay (minute) | Computational time (second) |
|---|---|---|---|---|---|
| Block 46 | 690-870 | 5 | 13590 | 215 | 3025 |
| Block 46's Tehran to Qom link | 690-870 | 5 | 13410 | 35 | 2898 |
| Block 45's Mashhad to Tehran link | 690-870 | 5 | 14510 | 1135 | 3206 |
| Block 451[1] | 690-870 | 5 | 14900 | 1525 | 3420 |
| Block 34 | 1450-1630 | 5 | 13655 | 280 | 2372 |
| Block 12 | 910-1090 | 5 | 13645 | 270 | 2722 |
| Block 23 | 1270-1450 | 5 | 14100 | 725 | 2347 |

The aforementioned results demonstrate that the impact of disruption duration on total delay varies depending on the disruption scenario and block's travel time. For instance, altering the disruption duration for Block 46's Tehran to Qom link almost has no effect on the total delay, as Train 1 can easily use its opposing block to complete its journey without significantly affecting Train 5's initial timetable.

The model's run time by the CPLEX solver increased significantly as the disruption duration increased, and it was so high for all of the above scenarios that it actually limited the usefulness of the model in the disruption management of a rail network, since the traffic flow in the network should come to a halt for a considerable amount of time before receiving an optimal post-disruption schedule. Therefore, an adjustment is proposed and tested in the second case study to encounter the run time issue.

## 6. Case Study 2: A Crowded Test Network

The traffic flow in the Iranian railway is moderate at its highest, as could be seen in the basic disruption scenario where only two trains were affected by the disruption. Hence, the second case study introduces a crowded test network to fully test the rescheduling optimization model. This network is depicted in Figure 7 and comprises 10 nodes and 10 blocks which are used by 25 trains to complete their trip in a 2-hour time horizon. Such a scenario could imply the operation of an urban rail transit during the daily peak period. The first and second disruption scenarios in this case study include the total closure of the two tracks of Block 67 and one track of Block 67, respectively.

---

[1] This scenario was infeasible initially because of the delay threshold parameter (β). The mentioned results in Table 8 were achieved after relaxing the corresponding constraint.

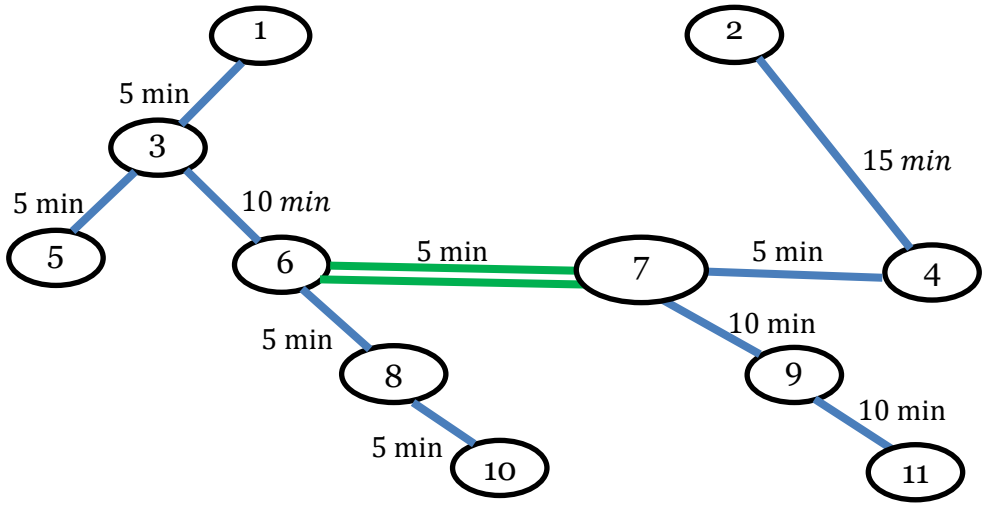

**Figure 7-** The test network under study in case study 2

Moreover, the model's parameters and the trains' initial timetable are presented in Tables 10 and 11.

**Table10**

The optimization model's parameters for the test network:

| Parameter | Value |
|---|---|
| $h_{ij}$ | 5 min |
| $DT_i$ | 0 min |
| $CAP_{ij}$ | 3 |
| $\delta$ | 5 min |
| $T$ | 120 min |
| $t_0$ | 0th min |
| $t_{dis}$ | 15th min |
| $\tau_{dis}$ | 20 min |
| $dis_{ij}^{gg'}$ (Disruption Scenario 1) | $\begin{cases} 2; \text{if:} ij = \{67\} \text{ and } gg' \in [15,35] \\ 0; \text{otherwise} \end{cases}$ |
| $dis_{ij}^{gg'}$ (Disruption Scenario 2) | $\begin{cases} 1; \text{if:} ij = \{67\} \text{ and } gg' \in [15,35] \\ 0; \text{otherwise} \end{cases}$ |
| $\beta$ | 15 min |

**Table11**

The initial timetable of trains in the test network:

| Train No. | Departure Time | Arrival Time | Total Travel Time | Origin | Destination |
|---|---|---|---|---|---|
| 1 | 0:00 | 0:10 | 10 | 1 | 5 |
| 2 | 0:00 | 0:10 | 10 | 7 | 9 |

| 3 | 0:00 | 0:20 | 20 | 7 | 5 |
|---|------|------|----|----|----|
| 4 | 0:10 | 0:25 | 15 | 10 | 7 |
| 5 | 0:10 | 0:25 | 15 | 4 | 2 |
| 6 | 0:10 | 0:20 | 10 | 6 | 4 |
| 7 | 0:10 | 0:30 | 20 | 1 | 8 |
| 8 | 0:10 | 0:35 | 25 | 9 | 10 |
| 9 | 0:15 | 0:55 | 40 | 11 | 2 |
| 10 | 0:20 | 0:40 | 20 | 4 | 10 |
| 11 | 0:20 | 0:40 | 20 | 10 | 3 |
| 12 | 0:25 | 0:35 | 10 | 5 | 1 |
| 13 | 0:25 | 1:05 | 40 | 2 | 5 |
| 14 | 0:20 | 0:55 | 35 | 10 | 11 |
| 15 | 0:30 | 0:50 | 20 | 2 | 7 |
| 16 | 0:30 | 0:45 | 15 | 6 | 1 |
| 17 | 0:35 | 1:15 | 40 | 1 | 2 |
| 18 | 0:35 | 0:55 | 20 | 11 | 7 |
| 19 | 1:00 | 1:20 | 20 | 3 | 10 |
| 20 | 0:50 | 1:00 | 10 | 1 | 5 |
| 21 | 0:55 | 1:15 | 20 | 11 | 7 |
| 22 | 1:00 | 1:15 | 15 | 10 | 7 |
| 23 | 1:05 | 1:35 | 30 | 5 | 9 |
| 24 | 1:10 | 1:30 | 20 | 4 | 10 |
| 25 | 1:15 | 1:55 | 40 | 2 | 5 |

Figure 8 is the space-time trajectory for the initial timetable, and the model's results for the disruption scenario 1 are depicted in Figure 9. The details of the results are mentioned in Table 12.

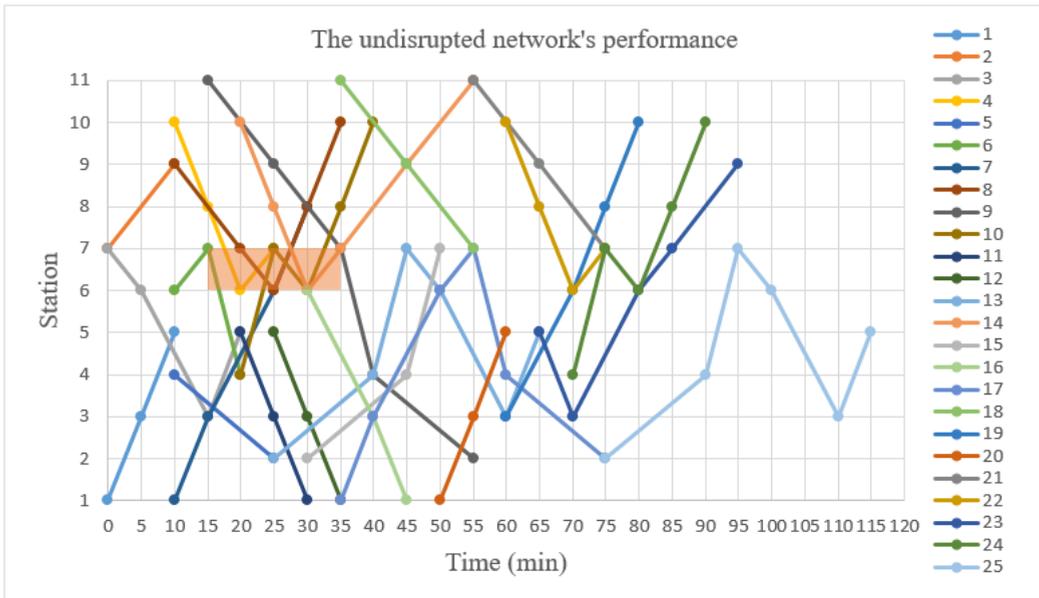

**Figure 8-** The space-time trajectory for the trains' initial timetable in the test network

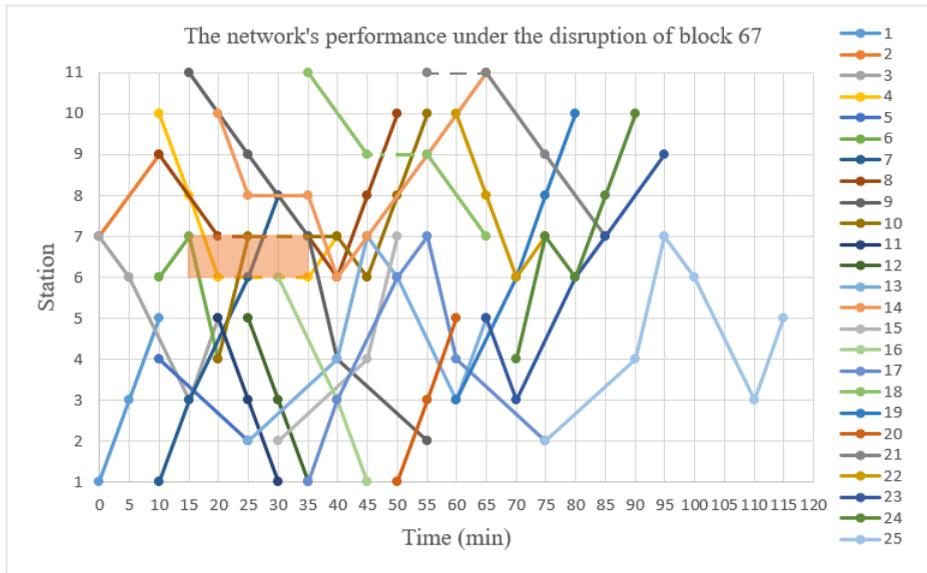

**Figure 9-** The space-time trajectory for the trains' timetable under disruption

**Table12**

The model's results for the test network under disruption scenarios 1 and 2:

| Disrupted block | Disruption period (minute) | $R^{*1}$ | The optimization model's objective function | Total delay (minute) | Computational time (second) |
|---|---|---|---|---|---|
| Block 67 | 15-35 | {4, 8, 10, 14} | 615 | 75 | 1593 |
| One of the Block 67's tracks | 15-35 | {4, 14} | 565 | 25 | 1282 |

As is mentioned in Table 12, the run time for both of the scenarios exceeds the 20-minute disruption period, which in turn, negates the effectiveness of the rescheduling model. Therefore, an adjusted model is introduced in the following section to decrease the run time.

## 7. Adjusted Rescheduling Model

In the adjusted model, the optimization model only focuses on the set of trains directly affected by the disruption ($R^*$), and the travel schedule of unaffected trains would not be altered. In this way, the complexity of the problem and the model's run time would decrease severely. To do so, a change should be made to the Constraints related to divergence from the given initial timetable (Equations 20-22). This change includes omitting Equation 20 from the model's constraints and replacing Equations 21 and 22 with 29 and 30, respectively.

---

[1] The set of trains directly affected by the disruption is created using the $\beta$ parameter. This set is comprised of the trains that their initial travel pattern makes them traverse the disrupted block during the disruption period and all of the trains that follow them into that block without having a buffer time of length $\beta$.

$$d_i^k - \bar{d}_i^k = 0; \; \forall i \in N^k \quad \forall t, t' \in T \quad \forall k \in \{R^+ \cup R^-\} - R^* \tag{29}$$

$$x_{ij}^{k,tt'} - \bar{x}_{ij}^{k,tt'} = 0; \; \forall ij \in A \quad \forall t, t' \in T \quad \forall k \in \{R^+ \cup R^-\} - R^* \tag{30}$$

The space-time trajectory resulting from the adjusted model for the previous disruption scenarios is depicted in Figures 10 and 11. The red dotted line in Figure 11 indicates that train 4 has traversed Block 67 using the opposite track. Moreover, the detailed results are presented in Table 13. As it is mentioned in Table 13, the runtime of the adjusted model has been reduced by up to 88% compared to the basic model; hence, the adjusted model can be easily implemented in real-world situations. Furthermore, the adjusted model may result in a better or worse outcome in comparison with the basic model, since the feasible region of the adjusted model is different from that of the basic model. In other words, the adjusted model is developed by relaxing a constraint (Equations 20, 21, and 22) and imposing other constraints (Equations 29 and 30) on the basic model.

**Table13**
The model's results for the test network under disruption scenarios 1 and 2:

| Disrupted block | Disruption period (minute) | $R^*$ | The optimization model's objective function | Total delay (minute) | Computational time (second) |
|---|---|---|---|---|---|
| Block 67 | 15-35 | {4, 8, 10, 14} | 605 | 65 | 224 |
| One of the Block 67's tracks | 15-35 | {4, 14} | 570 | 30 | 153 |

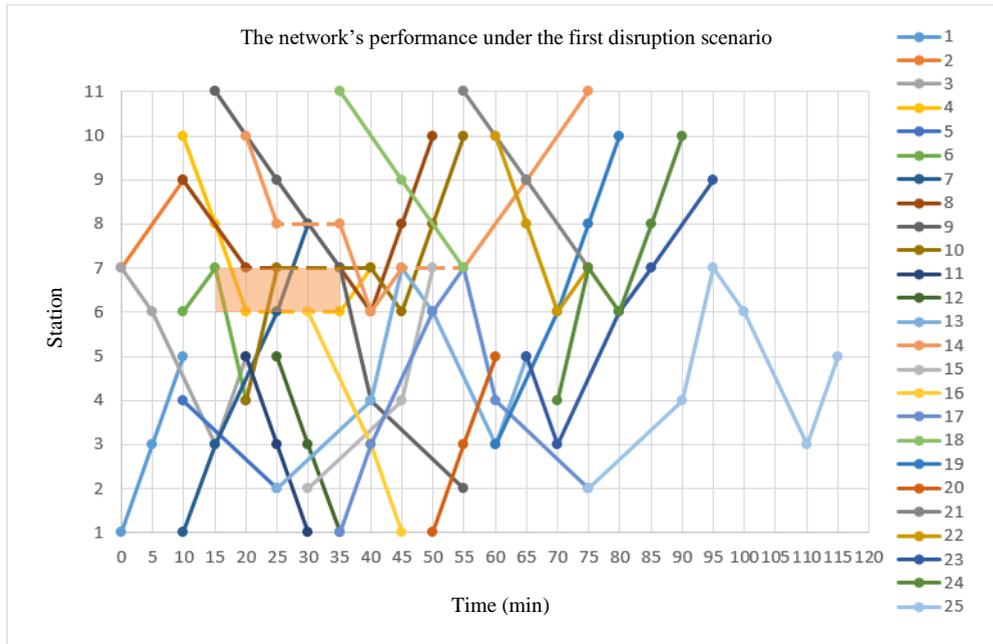

**Figure 10-** The space-time trajectory resulted from the adjusted model under the first disruption scenario

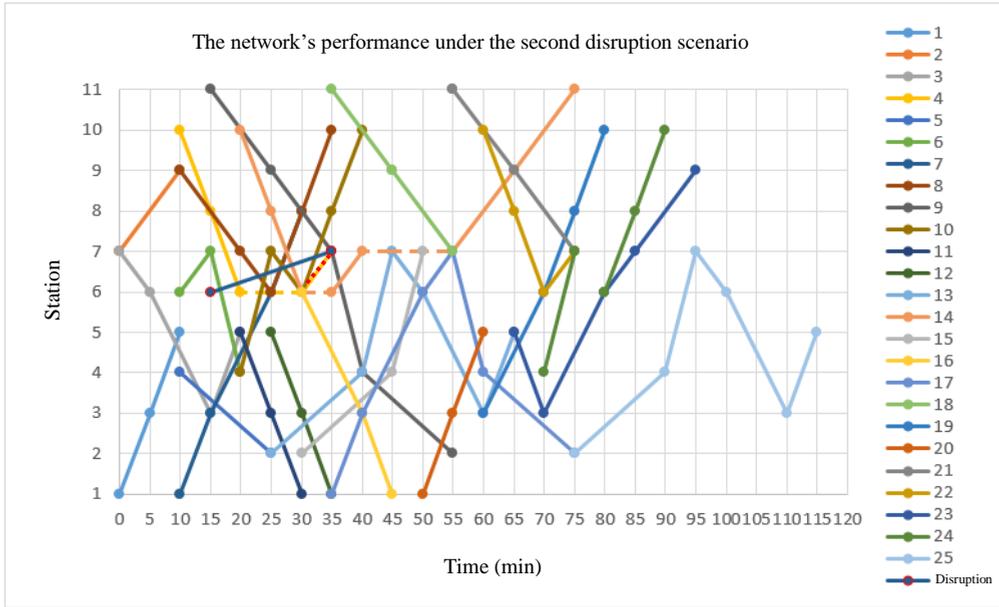

**Figure 11-** The space-time trajectory resulted from the adjusted model under the second disruption scenario

### 8. Conclusions

This study presents a novel methodology to mitigate disruptions in rail networks, which play a critical role in advancing digital infrastructure. The methodology provides an effective approach to address the rescheduling problem in railway networks by summarizing the network into a medium-scale representation, identifying its important nodes using both topological and passenger flow factors, and applying an optimization model for rescheduling on the summarized network. Two case studies were conducted, involving several disruption scenarios in a real-world network and a crowded test network. The results demonstrate that the model performs well in both sparse (Iranian railway) and crowded (test network) scenarios, achieving computational time reductions of up to 88% when adjustments to the model are implemented. This makes the approach highly applicable to real-world operations. Furthermore, the practical implications of implementing such a digitalized system include increased operational efficiency, reduced delays, benefits for passengers and freight operators, and enhanced system reliability and resilience.

Although the results of the case studies were plausible, several suggestions can enhance the model's performance and practicality. These include incorporating trip cancellations, developing heuristic solutions to address multiple disruptions simultaneously, and accounting for varying train velocities. While the model is capable of considering scenarios with multiple simultaneous disruptions, its runtime using the CPLEX solver is high. Therefore, future work should focus on developing heuristic solutions to mitigate this issue.

Additionally, the lack of consideration for freight trains in the model may affect the results, as freight train operators often have different objectives compared to passenger trains, affecting the available capacity for both. Future work could also address the potential negative impact of omitting demand from unimportant nodes, especially in dense networks where these nodes provide valuable opportunities for overtaking. This could be achieved by developing an alternative metric for the important nodes algorithm's cutoff.

Furthermore, rerouting processes and track adjustment strategies were not included in this work, as the model was designed for implementation on the summarized network generated by the important nodes algorithm, which is typically sparse with limited rerouting options. Finally, our methodology could be integrated with digital platforms for real-time decision-making in disrupted rail systems, provided the infrastructure is equipped with real-time data collection tools such as sensors for tracking train movements and automated fare collection systems. We recommend that future research explore this aspect to further enhance the model's applicability and effectiveness.

### Ethical Statement

The authors declare that they have no known competing financial interests or personal relationships that could have appeared to influence the work reported in this paper.

### References

[1] Poliński, J., & Ochociński, K. (2020). Digitization in rail transport. Problemy Kolejnictwa, (188), 137-148.

[2] Zhuravleva, N., Kliestik, T. (2023). Railway Transport Digitalization: Development Methodology and Effects of Digital Implementation Processes. In: Bencsik, A., Kulachinskaya, A. (eds) Digital Transformation: What is the Company of Today?. Lecture Notes in Networks and Systems, vol 805. Springer, Cham.

[3] Chechenova, L., & Batalova, N. (2023). Digitalization of freight rail transportation as a factor in improving their efficiency. E3S Web of Conferences, 460, 06042.

[4] Montero, J. (2020). The digitalization dilemma in the railway industry. In Handbook on Railway Regulation (pp. 379-396). Edward Elgar Publishing.

[5] W. Jing, X. Xu, Y. Pu, Route redundancy-based approach to identify the critical stations in metro networks: A mean-excess probability measure, Reliability Engineering & System Safety 204 (2020)107.

[6] L. Wang, S. Zheng, Y. Wang, L. Wang, Identification of critical nodes in multimodal transportation network, Physica A: Statistical Mechanics and its Applications 580 (2021) 126170.

[7] L. Sun, Y. Huang, Y. Chen, L. Yao, Vulnerability assessment of urban rail transit based on multi-static weighted method in Beijing, China, Transportation Research Part A: Policy and Practice 108 (2018) 12–24.

[8] N. Bešinović, R. F. Nassar, C. Szymula, Resilience assessment of railway networks: Combining infrastructure restoration and transport management, Reliability Engineering & System Safety 224 (2022) 108538.

[9] J. Tang, L. Xu, C. Luo, T. S. A. Ng, Multi-disruption resilience assessment of rail transit systems with optimized commuter flows, Reliability Engineering & System Safety 214 (2021) 107715.

[10] L. Xu, T. S. A. Ng, A robust mixed-integer linear programming model for mitigating rail transit disruptions under uncertainty, Transportation Science 54 (5) (2020) 1388–1407.

[11] L. Yang, X. Zhou, Z. Gao, Credibility-based rescheduling model in a double-track railway network: a fuzzy reliable optimization approach, Omega 48 (2014) 75–93.

[12] C. Gong, J. Shi, Y. Wang, H. Zhou, L. Yang, D. Chen, H. Pan, Train timetabling with dynamic and random passenger demand: A stochastic optimization method, Transportation Research Part C: Emerging Technologies 123 (2021) 102963.

[13] A. Bärmann, P. Gemander, L. Hager, F. Nöth, O. Schneider, EETTlib—energy-

efficient train timetabling library, Networks 81 (1) (2023) 51–74.

[14] P. G. Sessa, V. De Martinis, A. Bomhauer-Beins, U. A. Weidmann, F. Corman, A hybrid stochastic approach for offline train trajectory reconstruction, Public Transport 13 (2021) 675–698.

[15] Y. Wang, U. A. Weidmann, H. Wang, Using catastrophe theory to describe railway system safety and discuss system risk concept, Safety Science 91 (2017) 269–285.

[16] N. Leng, U. Weidmann, Discussions of the reschedule process of passengers, train operators and infrastructure managers in railway disruptions, Transportation Research Procedia 27 (2017) 538–544.

[17] S. Herrigel, M. Laumanns, J. Szabo, U. Weidmann, Periodic railway timetabling with sequential decomposition in the PESP model, Journal of rail transport planning & management 8 (3-4) (2018) 167–183.

[18] A. Toletti, M. Laumanns, U. Weidmann, Coordinated railway traffic rescheduling with the resource conflict graph model, Journal of Rail Transport Planning & Management 15 (2020) 100173.

[19] G. C. Caimi, Algorithmic decision support for train scheduling in a large and highly utilised railway network, ETH Zurich, (2009).

[20] N. Ghaemi, O. Cats, RMP. Goverde, Railway disruption management challenges and possible solution directions, Public Transport 9 (2017) 343–364.

[21] P. Kecman, F. Corman, A. D'Ariano, Rescheduling models for railway traffic management in large-scale networks, Public Transport 5 (2013) 95–123.

[22] RMP. Goverde, Railway timetable stability analysis using max-plus system theory, Transportation Research Part B: Methodological 41(2) (2007) 179–201.

[23] RMP. Goverde, A delay propagation algorithm for large-scale railway traffic networks, Transportation Research Part C, Emerging Technologies 18(3) (2010) 269–287.

[24] N. Bešinović, Resilience in railway transport systems: a literature review and research agenda, Transport Reviews 40 (4) (2020) 457–478.

[25] L. Ge, S. Voß, L. Xie, Robustness and disturbances in public transport. Public Transport 14 (2022) 191–261.